\documentclass[12pt]{article}
\usepackage{graphics,amsmath,amssymb}
\usepackage{latexsym}
\usepackage{epsfig}
\usepackage{hyperref}
\usepackage{fancybox}
\usepackage{mathtools}

\def\epsilon{\varepsilon}

\def\phi{\varphi}

\newtheorem{theorem}{Theorem}[section]

\def\R{{\mathbb R}}

\title{\bf Time fractional
 exact controllability}

\author{Paola Loreti
\thanks{Dipartimento di Scienze di Base e Applicate per l'Ingegneria,
Sapienza Universit\`a di Roma,
Via Antonio Scarpa 16, 00161 Roma (Italy); e-mail: 
$<$paola.loreti@uniroma1.it$>$ }
\and Daniela Sforza
\thanks{Dipartimento di Scienze di Base e Applicate per l'Ingegneria, 
Sapienza Universit\`a di Roma,
Via Antonio Scarpa 16, 00161 Roma (Italy); e-mail: 
$<$daniela.sforza@uniroma1.it$>$ }}

\begin{document}
\date{}

\maketitle

\begin{abstract}
Our purpose is to adapt the Hilbert Uniqueness Method by J.-L. Lions in the case of fractional diffusion-wave equations. The main difficulty is to determine the right shape for the adjoint system, suitable for the procedure of HUM.

\end{abstract}



\bigskip
\noindent
\section{Introduction}

In this paper we investigate the fractional differential equation
\begin{equation}\label{eq:problem-uI}
\sideset{^{C}}{_{0+}^{\alpha}}{\mathop D}u
=\Delta u
\end{equation}
where the symbol $\sideset{^{C}}{_{0+}^{\alpha}}{\mathop D}u$ denotes the Caputo derivative 
\begin{equation*}
\sideset{^{C}}{_{0+}^{\alpha}}{\mathop D}u(t)
=\frac1{\Gamma(2-\alpha)}\int_0^t (t-s)^{1-\alpha}u{''}(s)\ ds
\end{equation*}
defined in \cite{C}. We assume that the order $\alpha$ belongs to $(1,2)$, so \eqref{eq:problem-uI} is commonly called fractional diffusion-wave equation.
For a comparison between Caputo derivatives and Riemann-Liouville derivatives see the survey \cite{MG}.

In  \cite{SY} the well-posedness has been established by means of a Fourier series approach that involves Mittag-Leffler functions. 
Also in \cite{LoretiSforza2} the proof of trace regularity properties for the weak solutions of \eqref{eq:problem-uI} is based on the representation like Fourier expansions.
The paper \cite{M}  drew attention to the fact that  fractional wave equations rule the propagation of mechanical diffusive waves in viscoelastic media.

The aim of this paper is to illustrate the Hilbert Uniqueness Method for  \eqref{eq:problem-uI},
introduced by J.-L. Lions in \cite{Lio3}. Several scientists  show interest for HUM since it gives a general procedure to get exact controllability for distributed systems when the control acts on the boun\-da\-ry of the domain. 
This method is based on the idea that the observability of the adjoint problem is sufficient for the exact controllability of the assigned problem. The observability is also a necessary condition,  for detailed argumentations see \cite{K}.

Our task is to consider the following
problem: for a domain $\Omega$ and $T>0$, we have to control the system
\begin{equation}\label{eq:phiI}
\left \{\begin{array}{l}\displaystyle
\sideset{^{C}}{_{0+}^{\alpha}}{\mathop D}u
=\Delta u
\qquad \text{in }\  (0,T)\times\Omega,
\\
\\
u(0)=u_{t}(0)=0\qquad  \text{in }\ \Omega\,,
\end{array}\right .
\end{equation} 
by means of a function that acts on a  boundary's subset.

We point out that, due to the nonlocal nature of the problem,  HUM allows one to achieve  exact reachability rather than exact controllability.

The main difficulty in carrying out HUM for  \eqref{eq:phiI}   is to find the right form for the adjoint system. Precisely, we have to change the Caputo derivative $\sideset{^{C}}{_{0+}^{\alpha}}{\mathop D}$ into the Riemann-Liouville derivative $\sideset{}{_{T-}^{\alpha}}{\mathop D}$ and consider as adjoint system associated with \eqref{eq:phiI} the following problem
\begin{equation*}
\begin{cases}
\displaystyle 
\sideset{}{_{T-}^{\alpha}}{\mathop D}w
=\Delta w
\qquad \text{in }\  (0,T)\times\Omega,
\\
\\
w=0, 
\qquad \text{on }\  (0,T)\times\partial\Omega,
\end{cases}
\end{equation*}
with  final data 
\begin{equation*} 
\sideset{}{_{T-}^{\alpha-2}}{\mathop D}w(T)=w_{0}\,,\qquad \sideset{}{_{T-}^{\alpha-1}}{\mathop D}w(T)=w_{1}
\qquad  \text{in }\ \Omega\,.
\end{equation*} 
We remark that the question of the exact reachability is still unanswered, since inverse and direct inequalities have to be proved.

Another adaptation of HUM to nonlocal problems with regular integral kernels can be found in \cite{LoretiSforza1}.

In Section 2 we set up notations and terminology. Section 3 provides a detailed exposition of the Hilbert Uniqueness Method for system \eqref{eq:phiI}.

\section{Preliminaries}
Let $f\in L^1(0,T)$ $(T>0)$ be.
The following definitions for fractional integrals and derivatives will be used throughout the paper. 
\begin{itemize}
\item[{\bf --}] 
{\it Riemann--Liouville fractional integrals} $I_{0+}^{\alpha}f$ and $I_{T-}^{\alpha}f$ of order $\alpha>0$:
\begin{equation}
I_{0+}^{\alpha}f(t)=\frac1{\Gamma(\alpha)}\int_0^t (t-s)^{\alpha-1}f(s)\ ds
\qquad \mbox{a.e.}\ t\in(0,T),
\end{equation}
\begin{equation}
I_{T-}^{\alpha}f(t)=\frac1{\Gamma(\alpha)}\int_t^T (s-t)^{\alpha-1}f(s)\ ds,
\qquad \mbox{a.e.}\ t\in(0,T),
\end{equation}
where $\Gamma (\alpha)=\int_0^\infty t^{\alpha-1}e^{-t}\ dt$ is the Euler gamma function. 
\item[{\bf --}] 
{\it Caputo fractional derivative} $\sideset{^{C}}{_{0+}^{\alpha}}{\mathop D}f$
of order $\alpha\in (1,2)$ when $f'$ is absolutely continuous: 
\begin{equation}
\sideset{^{C}}{_{0+}^{\alpha}}{\mathop D}f(t)
=I_{0+}^{2-\alpha}f{''}(t)
=\frac1{\Gamma(2-\alpha)}\int_0^t (t-s)^{1-\alpha}f{''}(s)\ ds
\,.
\end{equation}
\item[{\bf --}] 
{\it Riemann--Liouville fractional derivative} 
$\sideset{}{_{T-}^{\alpha}}{\mathop D}f$ of order $\alpha\in (n-1,n)$,  $n=0,1,2$:
\begin{equation}\label{eq:RL-T}
\sideset{}{_{T-}^{\alpha}}{\mathop D}f(t)=\Big(\frac{d}{dt}\Big)^nI_{T-}^{n-\alpha}f(t)
=\frac1{\Gamma(n-\alpha)}\Big(\frac{d}{dt}\Big)^n\int_t^T (s-t)^{n-\alpha-1}f(s)\ ds.
\end{equation}
\end{itemize}
In particular, we note that for $-1<\alpha<0$, that is $n=0$, we have
\begin{equation}\label{eq:n=0}
\sideset{}{_{T-}^{\alpha}}{\mathop D}f=I_{T-}^{-\alpha}f.
\end{equation}
We will use later a well-known property (see e.g. \cite{Kilbas2006}):
\begin{equation}\label{eq:change+-}
\int_0^T I_{0+}^{\alpha}f(t)g(t)\ dt=\int_0^T f(t)I_{T-}^{\alpha}g(t)\ dt\,.
\end{equation}

\section{Time fractional exact controllability}\label{se:HUM}
Our aim is to adapt the Hilbert Uniqueness Method by J.-L. Lions \cite{Lio3} for fractional differential equations.

Let $T>0$ and $\Omega$ be a bounded domain of class $C^2$ in $\R^N$, $N\ge1$, with boundary $\partial\Omega$ given by 
$\partial\Omega=\partial\Omega_1\cup\partial\Omega_2$ and $\partial\Omega_1\cap\partial\Omega_2=\emptyset$.
We consider the fractional diffusion-wave equation with $\alpha\in(1,2)$:
\begin{equation}\label{eq:problem-u}
\sideset{^{C}}{_{0+}^{\alpha}}{\mathop D}u(t,x)
=\Delta u(t,x)
\qquad
(t,x)\in (0,T)\times\Omega,
\end{equation}
with null initial conditions 
\begin{equation}
u(0,x)=u_{t}(0,x)=0\qquad  x\in\Omega\,,
\end{equation} 
and boundary conditions, 
\begin{equation}\label{eq:bound-u1}
u(t,x)=g(t,x) \quad (t,x)\in (0,T)\times\partial\Omega_1,
\qquad u(t,x)=0
\quad (t,x)\in (0,T)\times\partial\Omega_2\,.
\end{equation}
We give the definition for a reachability problem: fixed
$u_{0}\in  L^{2}(\Omega)$ and $u_{1}\in   H^{-1}(\Omega)$,
we have to determine $g\in L^2(0,T;L^2(\partial\Omega_1) )$ such that the weak
solution $u$ of problem  \eqref{eq:problem-u}-\eqref{eq:bound-u1}
verifies the final conditions
\begin{equation}\label{eq:problem-u1}
u(T,x)=u_{0}(x)\,,\quad u_{t}(T,x)=u_{1}(x)\,,
\qquad x\in\Omega
\,.
\end{equation}
One can solve such type of problems by the Hilbert Uniqueness Method. To this end, we proceed as follows.
Given
$(w_{0},w_{1})\in (C^\infty_c(\Omega))^2$,
we introduce the {\it adjoint} system of (\ref{eq:problem-u}), that is 
\begin{equation}\label{eq:adjoint}
\begin{cases}
\displaystyle 
\sideset{}{_{T-}^{\alpha}}{\mathop D}w(t,x)
=\Delta w(t,x)
\qquad
(t,x)\in (0,T)\times\Omega,
\\
\\
w(t,x)=0 
\qquad (t,x)\in (0,T)\times\partial\Omega
\,,
\end{cases}
\end{equation}
with  final data 
\begin{equation} \label{eq:final}
\sideset{}{_{T-}^{\alpha-2}}{\mathop D}w(T,\cdot)=w_{0}\,,\qquad \sideset{}{_{T-}^{\alpha-1}}{\mathop D}w(T,\cdot)=w_{1}\,.
\end{equation} 
The above problem is well-posed by means of a spectral representation and the Mittag-Leffler functions, see e.g. \cite{Kilbas2006,HP,IK,LPJ,HAR}. 

Thanks to the regularity of the final data, the solution  $w$ of \eqref{eq:adjoint}--\eqref{eq:final}
is regular enough to 
take in \eqref{eq:bound-u1} $g=\partial_\nu w$. So  the  nonhomogeneous problem \eqref{eq:problem-u}-\eqref{eq:bound-u1} can be written in the form
\begin{equation}\label{eq:phi}
\left \{\begin{array}{l}\displaystyle
\sideset{^{C}}{_{0+}^{\alpha}}{\mathop D}u(t,x)
=\Delta u(t,x)
\qquad (t,x)\in (0,T)\times\Omega,
\\
\\
u(0,x)=u_{t}(0,x)=0\qquad  x\in\Omega\,,
\\
\displaystyle
u(t,x)=\partial_\nu w(t,x), \ \  (t,x)\in (0,T)\times\partial\Omega_1,
\quad 
u(t,x)=0
\ \ (t,x)\in (0,T)\times\partial\Omega_2\,.
\end{array}\right .
\end{equation} 
As in the non integral case, it can be proved  that  problem \eqref{eq:phi} admits a unique solution
$u$.
We can introduce a linear operator defined as 
\begin{equation}\label{eq:psi0}
\Psi(w_{0},w_{1})=\big(-u_t(T,\cdot),u(T,\cdot)\big),
\qquad
(w_{0},w_{1})\in (C^\infty_c(\Omega))^2
\,.
\end{equation}
We will prove that
\begin{equation}\label{eq:psi}
\langle\Psi(w_{0},w_{1}),(w_{0},w_{1})\rangle_{L^2(\Omega)}
=\int_0^T\int_{\partial\Omega_1}\big|\partial_{\nu}w(t,x)\big|^2 \ dx dt
 \,. 
\end{equation}
To this end, we multiply the  equation in (\ref{eq:phi}) by the solution $w(t,x)$ of the adjoint system \eqref{eq:adjoint}-\eqref{eq:final} and integrate 
on $(0,T)\times\Omega$, that is
\begin{equation}\label{eq:mult}
\int_0^T\int_{\Omega}\sideset{^{C}}{_{0+}^{\alpha}}{\mathop D}u(t,x)w(t,x)\ dx dt
=\int_0^T\int_{\Omega}\Delta u(t,x)w(t,x)\ dx dt
\,.
\end{equation}
If we apply \eqref{eq:change+-}, then we obtain for any $x\in\Omega$
\begin{equation}\label{eq:intp}
\int_0^T \sideset{^{C}}{_{0+}^{\alpha}}{\mathop D}u(t,x) w(t,x)\ dt
=\int_0^TI_{0+}^{2-\alpha}u_{tt}(t,x) w(t,x)\ dt
=\int_0^Tu_{tt}(t,x) I_{T-}^{2-\alpha}w(t,x)\ dt
\,.
\end{equation}
Integrating twice by parts and taking into account that $u(0,\cdot)=u_t(0,\cdot)=0$, we have
\begin{multline}\label{eq:intp1}
\int_0^Tu_{tt}(t,x)I_{T-}^{2-\alpha}w(t,x)dt
=u_t(T,x)I_{T-}^{2-\alpha} w(T,x)-\int_0^Tu_t(t,x)\frac{\partial}{\partial t}I_{T-}^{2-\alpha}w(t,x)dt
\\
=u_t(T,x)I_{T-}^{2-\alpha} w(T,x)-u(T,x)\frac{\partial}{\partial t}I_{T-}^{2-\alpha}w(T,x)
+\int_0^Tu(t,x)\sideset{}{_{T-}^{\alpha}}{\mathop D}w(t,x)   \ dt
\,.
\end{multline}
Since $-1<\alpha-2<0$, we can use \eqref{eq:n=0} to get
\begin{equation*}
I_{T-}^{2-\alpha}w(T,x)
=\sideset{}{_{T-}^{\alpha-2}}{\mathop D}w(T,x).
\end{equation*}
Moreover, $0<\alpha-1<1$, so by \eqref{eq:RL-T} for $n=1$ we have
\begin{equation*}
\frac{\partial}{\partial t}I_{T-}^{2-\alpha}w(T,x)
=\sideset{}{_{T-}^{\alpha-1}}{\mathop D}w(T,x)
\,.
\end{equation*}
Combining \eqref{eq:intp} with \eqref{eq:intp1} yields
\begin{multline}\label{eq:change1}
\int_0^T\sideset{^{C}}{_{0+}^{\alpha}}{\mathop D}u(t,x)w(t,x)\ dt
\\
=u_t(T,x)\sideset{}{_{T-}^{\alpha-2}}{\mathop D}w(T,x)-u(T,x)\sideset{}{_{T-}^{\alpha-1}}{\mathop D}w(T,x)
+\int_0^Tu(t,x)\sideset{}{_{T-}^{\alpha}}{\mathop D}w(t,x)   \ dt
\,.
\end{multline}
On the other hand, integrating twice by parts with respect to the variable $x$ we have
\begin{equation}\label{eq:change2}
\int_{\Omega}\Delta u(t,x)w(t,x) dx=
\int_\Omega u(t,x)\triangle w(t,x)\ dx
-\int_{\partial\Omega}u(t,x)\partial_{\nu}w(t,x) \ dx
 \,.
\end{equation}
Putting \eqref{eq:change1}  and  \eqref{eq:change2} into \eqref{eq:mult} we get
\begin{multline*}
\int_0^T\int_\Omega u(t,x)\sideset{}{_{T-}^{\alpha}}{\mathop D}w(t,x) \ dx dt 
+\int_\Omega u_t(T,x)\sideset{}{_{T-}^{\alpha-2}}{\mathop D} w(T,x) dx
-\int_\Omega u(T,x)\sideset{}{_{T-}^{\alpha-1}}{\mathop D}w(T,x) dx
\\
=\int_0^T\int_\Omega u(t,x)\triangle w(t,x) \ dx dt
-\int_0^T\int_{\partial\Omega}u(t,x)\partial_{\nu}w(t,x) \ dx dt
\,.
\end{multline*}
We recall that $w$ is the solution of the adjont problem \eqref{eq:adjoint}:
keeping  
$
\sideset{}{_{T-}^{\alpha}}{\mathop D}w
=\Delta w
$
in mind, 
from the above equation we have
\begin{equation*}
\int_\Omega u_t(T,x)\sideset{}{_{T-}^{\alpha-2}}{\mathop D} w(T,x) dx
-\int_\Omega u(T,x)\sideset{}{_{T-}^{\alpha-1}}{\mathop D}w(T,x) dx
=-\int_0^T\int_{\partial\Omega}u(t,x)\partial_{\nu}w(t,x) \ dx dt
\,.
\end{equation*}
By the final data 
\eqref{eq:final} of $w$ and the boundary conditions satisfied by $u$, see \eqref{eq:phi}, we deduce that
\begin{multline*}
\langle\Psi(w_{0},w_{1}),(w_{0},w_{1})\rangle_{L^2(\Omega)}
=
\int_\Omega  u(T,x) w_1dx
-\int_\Omega  u_t(T,x)w_0 dx
=\int_0^T\int_{\partial\Omega_1}\big|\partial_{\nu}w(t,x)\big|^2 \ dx dt
\,,
\end{multline*}
hence \eqref{eq:psi}.

Define a semi-norm on the space $\big(C^\infty_c(\Omega)\big)^2$:
for $(w_{0},w_{1})\in \big(C^\infty_c(\Omega)\big)^2$ we set
\begin{equation}\label{eq:normF}
\|(w_{0},w_{1})\|^2_{F}:=
\int_0^T\int_{\partial\Omega_1}\big|\partial_{\nu}w(t,x)\big|^2 \ dx dt\,.
\end{equation}
We observe that 
$\|\cdot\|_{F}$ is a norm if and only if the following uniqueness theorem holds.
\begin{theorem}\label{th:uniqueness}
If $w$ is the solution of problem {\rm (\ref{eq:adjoint})--(\ref{eq:final})} such that
$$
\partial_{\nu}w(t,x)=0\,,\qquad \forall (t,x)\in [0,T]\times\partial\Omega_1\,,
$$
then 
$$
w(t,x)= 0 \qquad\forall (t,x)\in [0,T]\times\Omega\,.
$$
\end{theorem}
If theorem \ref{th:uniqueness} holds true, then we can define the Hilbert space $F$ as the completion of $ \big(C^\infty_c(\Omega)\big)^2$ for
the norm (\ref{eq:normF}). Moreover, the operator $\Psi$ extends uniquely to a continuous operator, denoted again by $\Psi$, from $F$ to the dual
space $F'$ in such a way that   $\Psi:F\to F'$ is an isomorphism.

In conclusion,
if we prove  the uniqueness result given by theorem \ref{th:uniqueness} and 
$$F=H^1_0(\Omega)\times L^2(\Omega)$$
with the equivalence of the respective norms, then we can solve the reachability problem
\eqref{eq:problem-u}--\eqref{eq:problem-u1}
taking $(u_{0},u_{1})\in  L^{2}(\Omega)\times H^{-1}(\Omega)$.


\begin{thebibliography}{99}

\itemsep=\smallskipamount



\bibitem{C} M. Caputo, Linear model of dissipation whose Q is almost frequency independent-II, in Geophysical Journal International, 13  (1967) 529--539. 

\bibitem{HP}  N. Heymans, I. Podlubny, 
Physical interpretation of initial conditions for fractional differential equations with Riemann-Liouville fractional derivatives. Rheologica Acta, 45 (2006), 765--771. 

\bibitem{HAR}
S. Hristova, R. Agarwal, D. O'Regan, Explicit solutions of initial value problems for systems of linear Riemann-Liouville fractional differential equations with constant delay. Adv. Difference Equ. (2020), Paper No. 180, 18 pp. 

\bibitem{IK} D. Idczak, R. Kamocki, On the existence and uniqueness and formula for the solution of R-L fractional Cauchy problem in $\R^n$. Fract. Calc. Appl. Anal. 14 (2011), 538--553. 


\bibitem{Kilbas2006} 
A.A. Kilbas, H.M. Srivastava, J.J. Trujillo, 
 \emph{Theory and Applications of Fractional Differential Equations}. 
 Elsevier, Amsterdam (2006).


\bibitem{K} V. Komornik, {\em Exact controllability and stabilization. The multiplier method}, RAM: Research in Applied Mathematics, Masson, Paris; John Wiley and Sons, Ltd., Chichester, 1994.


\bibitem{LPJ}
K. Li, J. Peng, J. Jia, 
Cauchy problems for fractional differential equations with Riemann-Liouville fractional derivatives. 
J. Funct. Anal. 263 (2012), 476--510. 



\bibitem{Lio3} J.-L. Lions,
{\em Contr\^olabilit\'e exacte, perturbations et stabilisation de syst\`emes distribu\'es I-II}, 
Recherches en Math\'ematiques Appliqu\'ees, 8-9 (1988), Masson, Paris.

\bibitem{LoretiSforza1} P. Loreti, D. Sforza,  {\em Reachability problems 
for a class of integro-differential equations},  J. Differential 
Equations   248 (2010), 1711--1755.

\bibitem{LoretiSforza2} P. Loreti, D. Sforza,  {\em Fractional diffusion-wave equations: hidden regularity for weak solutions.} Fract. Calc. Appl. Anal. 24 (2021), 1015--1034.


\bibitem{MG} F. Mainardi, R. Gorenflo,
Time-fractional derivatives in relaxation processes: a tutorial survey. 
Fract. Calc. Appl. Anal. 10 (2007), 269--308. 

\bibitem{M} F. Mainardi, Fractional diffusive waves in viscoelastic solids, in: J.L. Wegner, F.R. Norwood (Eds.), Nonlinear Waves in Solids, ASME/AMR, Fairfield (1995) 93--97.

\bibitem{SY}
K. Sakamoto, M. Yamamoto,
Initial value/boundary value problems for fractional diffusion-wave equations and applications to some inverse problems.
\emph{J. Math. Anal. Appl.} \textbf{382}, No 1 (2011), 426--447.
\end{thebibliography}
\end{document}